# Hamilton and long cycles in $t$-tough graphs with $t > 1$

Zh.G. Nikoghosyan*

February 29, 2012

**Abstract**

It is proved that if $G$ is a $t$-tough graph of order $n$ and minimum degree $\delta$ with $t > 1$ then either $G$ has a cycle of length at least $\min\{n, 2\delta + 4\}$ or $G$ is the Petersen graph.

## 1 Introduction

Only finite undirected graphs without loops or multiple edges are considered. We reserve $n$, $\delta$, $\kappa$, $c$ and $\tau$ to denote the number of vertices (order), the minimum degree, connectivity, circumference and the toughness of a graph, respectively. A good reference for any undefined terms is [2].

The earliest lower bound for the circumference was developed in 1952 due to Dirac [3].

**Theorem A [3].** In every 2-connected graph, $c \geq \min\{n, 2\delta\}$.

In 1986, Bauer and Schmeichel [1] proved that the bound $2\delta$ in Theorem A can be enlarged to $2\delta + 2$ by replacing the 2-connectivity condition with 1-toughness.

**Theorem B [1].** In every 1-tough graph, $c \geq \min\{n, 2\delta + 2\}$.

In this paper we prove that in Theorem B the bound $2\delta + 2$ itself can be enlarged up to $2\delta + 4$ if $\tau > 1$ and $G$ is not the Petersen graph.

**Theorem 1.** Let $G$ be a graph with $\tau > 1$. Then either $c \geq \min\{n, 2\delta + 4\}$ or $G$ is the Petersen graph.

---

*G.G. Nicoghossian (up to 1997)

To prove Theorem 1, we need the following result due to Voss [4].

**Theorem C [4]**. Let $G$ be a hamiltonian graph, $\{v_1, v_2, ..., v_t\} \subseteq V(G)$ and $d(v_i) \geq t$ $(i = 1, 2, ..., t)$. Then each pair $x, y$ of vertices of $G$ is connected in $G$ by a path of length at least $t$.

## 2 Notations and preliminaries

The set of vertices of a graph $G$ is denoted by $V(G)$ and the set of edges by $E(G)$. For $S$ a subset of $V(G)$, we denote by $G \backslash S$ the maximum subgraph of $G$ with vertex set $V(G) \backslash S$. We write $G[S]$ for the subgraph of $G$ induced by $S$. For a subgraph $H$ of $G$ we use $G \backslash H$ short for $G \backslash V(H)$. The neighborhood of a vertex $x \in V(G)$ will be denoted by $N(x)$. Furthermore, for a subgraph $H$ of $G$ and $x \in V(G)$, we define $N_H(x) = N(x) \cap V(H)$ and $d_H(x) = |N_H(x)|$. Let $s(G)$ denote the number of components of a graph $G$. A graph $G$ is $t$-tough if $|S| \geq ts(G \backslash S)$ for every subset $S$ of the vertex set $V(G)$ with $s(G \backslash S) > 1$. The toughness of $G$, denoted $\tau(G)$, is the maximum value of $t$ for which $G$ is $t$-tough (taking $\tau(K_n) = \infty$ for all $n \geq 1$).

A simple cycle (or just a cycle) $C$ of length $t$ is a sequence $v_1 v_2 ... v_t v_1$ of distinct vertices $v_1, ..., v_t$ with $v_i v_{i+1} \in E(G)$ for each $i \in \{1, ..., t\}$, where $v_{t+1} = v_1$. When $t = 2$, the cycle $C = v_1 v_2 v_1$ on two vertices $v_1, v_2$ coincides with the edge $v_1 v_2$, and when $t = 1$, the cycle $C = v_1$ coincides with the vertex $v_1$. So, all vertices and edges in a graph can be considered as cycles of lengths 1 and 2, respectively. A graph $G$ is hamiltonian if $G$ contains a Hamilton cycle, i.e. a cycle of length $n$. A cycle $C$ in $G$ is dominating if $G \backslash C$ is edgeless.

Paths and cycles in a graph $G$ are considered as subgraphs of $G$. If $Q$ is a path or a cycle, then the length of $Q$, denoted $|Q|$, is $|E(Q)|$. We write $Q$ with a given orientation by $\overrightarrow{Q}$. For $x, y \in V(Q)$, we denote by $x \overrightarrow{Q} y$ the subpath of $Q$ in the chosen direction from $x$ to $y$. For $x \in V(C)$, we denote the $h$-th successor and the $h$-th predecessor of $x$ on $\overrightarrow{C}$ by $x^{+h}$ and $x^{-h}$, respectively. We abbreviate $x^{+1}$ and $x^{-1}$ by $x^+$ and $x^-$, respectively. For each $X \subset V(C)$, we define $X^{+h} = \{x^{+h} | x \in X\}$ and $X^{-h} = \{x^{-h} | x \in X\}$.

**Special definitions**. Let $G$ be a graph, $C$ a longest cycle in $G$ and $P = x \overrightarrow{P} y$ a longest path in $G \backslash C$ of length $\overline{p} \geq 0$. Let $\xi_1, \xi_2, ..., \xi_s$ be the elements of $N_C(x) \cup N_C(y)$ occuring on $C$ in a consecutive order. Set

$$I_i = \xi_i \overrightarrow{C} \xi_{i+1}, \ I_i^* = \xi_i^+ \overrightarrow{C} \xi_{i+1}^- \ (i = 1, 2, ..., s),$$

where $\xi_{s+1} = \xi_1$.

(1) The segments $I_1, I_2, ..., I_s$ are called elementary segments on $C$ created by $N_C(x) \cup N_C(y)$.

(2) We call a path $L = z \overrightarrow{L} w$ an intermediate path between two distinct

elementary segments $I_a$ and $I_b$ if

$$z \in V(I_a^*), \ w \in V(I_b^*), \ V(L) \cap V(C \cup P) = \{z, w\}.$$

(3) Define $\Upsilon(I_{i_1}, I_{i_2}, ..., I_{i_t})$ to be the set of all intermediate paths between elementary segments $I_{i_1}, I_{i_2}, ..., I_{i_t}$.

**Lemma 1.** Let $G$ be a graph, $C$ a longest cycle in $G$ and $P = x\overrightarrow{P}y$ a longest path in $G \backslash C$ of length $\overline{p} \geq 1$. If $|N_C(x)| \geq 2$, $|N_C(y)| \geq 2$ and $N_C(x) \neq N_C(y)$ then

$$|C| \geq \begin{cases} 3\delta + \max\{\sigma_1, \sigma_2\} - 1 \geq 3\delta & \text{if} \quad \overline{p} = 1, \\ \max\{2\overline{p} + 8, 4\delta - 2\overline{p}\} & \text{if} \quad \overline{p} \geq 2, \end{cases}$$

where $\sigma_1 = |N_C(x) \backslash N_C(y)|$ and $\sigma_2 = |N_C(y) \backslash N_C(x)|$.

**Lemma 2.** Let $G$ be a graph, $C$ a longest cycle in $G$ and $P = x\overrightarrow{P}y$ a longest path in $G \backslash C$ of length $\overline{p} \geq 0$. If $N_C(x) = N_C(y)$ and $|N_C(x)| \geq 2$ then for each elementary segments $I_a$ and $I_b$ induced by $N_C(x) \cup N_C(y)$,

(a1) if $L$ is an intermediate path between $I_a$ and $I_b$ then

$$|I_a| + |I_b| \geq 2\overline{p} + 2|L| + 4,$$

(a2) if $\Upsilon(I_a, I_b) \subseteq E(G)$ and $|\Upsilon(I_a, I_b)| = i$ for some $i \in \{1, 2, 3\}$ then

$$|I_a| + |I_b| \geq 2\overline{p} + i + 5,$$

(a3) if $\Upsilon(I_a, I_b) \subseteq E(G)$ and $\Upsilon(I_a, I_b)$ contains two independent intermediate edges then

$$|I_a| + |I_b| \geq 2\overline{p} + 8.$$

**Lemma 3.** Let $G$ be a graph and $C$ a longest cycle in $G$. Then either $|C| \geq \kappa(\delta + 1)$ or there is a longest path $P = x_1\overrightarrow{P}x_2$ in $G \backslash C$ with $|N_C(x_i)| \geq 2$ ($i = 1, 2$).

## 3 Proofs

**Proof of Lemma 1.** Put

$$A_1 = N_C(x) \backslash N_C(y), \ A_2 = N_C(y) \backslash N_C(x), \ M = N_C(x) \cap N_C(y).$$

By the hypothesis, $N_C(x) \neq N_C(y)$, implying that

$$\max\{|A_1|, |A_2|\} \geq 1.$$

Let $\xi_1, \xi_2, ..., \xi_s$ be the elements of $N_C(x) \cup N_C(y)$ occuring on $C$ in a consecutive order. Put $I_i = \xi_i \overrightarrow{C} \xi_{i+1}$ ($i = 1, 2, ..., s$), where $\xi_{s+1} = \xi_1$. Clearly, $s = |A_1| + |A_2| + |M|$. Since $C$ is extreme, $|I_i| \geq 2$ ($i = 1, 2, ..., s$). Next, if $\{\xi_i, \xi_{i+1}\} \cap M \neq \emptyset$ for some $i \in \{1, 2, ..., s\}$ then $|I_i| \geq \overline{p} + 2$. Further, if either $\xi_i \in A_1, \xi_{i+1} \in A_2$ or $\xi_i \in A_2, \xi_{i+1} \in A_1$ then again $|I_i| \geq \overline{p} + 2$.

**Case 1**. $\overline{p} = 1$.
**Case 1.1**. $|A_i| \geq 1$ $(i = 1, 2)$.
It follows that among $I_1, I_2, ..., I_s$ there are $|M| + 2$ segments of length at least $\overline{p} + 2$. Observing also that each of the remaining $s - (|M| + 2)$ segments has a length at least 2, we have

$$|C| \geq (\overline{p} + 2)(|M| + 2) + 2(s - |M| - 2)$$

$$= 3(|M| + 2) + 2(|A_1| + |A_2| - 2)$$

$$= 2|A_1| + 2|A_2| + 3|M| + 2.$$

Since $|A_1| = d(x) - |M| - 1$ and $|A_2| = d(y) - |M| - 1$,

$$|C| \geq 2d(x) + 2d(y) - |M| - 2 \geq 3\delta + d(x) - |M| - 2.$$

Recalling that $d(x) = |M| + |A_1| + 1$, we get

$$|C| \geq 3\delta + |A_1| - 1 = 3\delta + \sigma_1 - 1.$$

Analogously, $|C| \geq 3\delta + \sigma_2 - 1$. So,

$$|C| \geq 3\delta + \max\{\sigma_1, \sigma_2\} - 1 \geq 3\delta.$$

**Case 1.2**. Either $|A_1| \geq 1, |A_2| = 0$ or $|A_1| = 0, |A_2| \geq 1$.
Assume w.l.o.g. that $|A_1| \geq 1$ and $|A_2| = 0$, i.e. $|N_C(y)| = |M| \geq 2$ and $s = |A_1| + |M|$. Hence, among $I_1, I_2, ..., I_s$ there are $|M| + 1$ segments of length at least $\overline{p} + 2 = 3$. Taking into account that each of the remaining $s - (|M| + 1)$ segments has a length at least 2 and $|M| + 1 = d(y)$, we get

$$|C| \geq 3(|M| + 1) + 2(s - |M| - 1) = 3d(y) + 2(|A_1| - 1)$$

$$\geq 3\delta + |A_1| - 1 = 3\delta + \max\{\sigma_1, \sigma_2\} - 1 \geq 3\delta.$$

**Case 2**. $\overline{p} \geq 2$.
We first prove that $|C| \geq 2\overline{p} + 8$. Since $|N_C(x)| \geq 2$ and $|N_C(y)| \geq 2$, there are at least two segments among $I_1, I_2, ..., I_s$ of length at least $\overline{p} + 2$. If $|M| = 0$ then clearly $s \geq 4$ and

$$|C| \geq 2(\overline{p} + 2) + 2(s - 2) \geq 2\overline{p} + 8.$$

Otherwise, since $\max\{|A_1|, |A_2|\} \geq 1$, there are at least three elementary segments of length at least $\overline{p} + 2$, that is

$$|C| \geq 3(\overline{p} + 2) \geq 2\overline{p} + 8.$$

So, in any case, $|C| \geq 2\overline{p} + 8$.

To prove that $|C| \geq 4\delta - 2\overline{p}$, we distinguish two main cases.

**Case 2.1**. $|A_i| \geq 1$ $(i = 1, 2)$.

It follows that among $I_1, I_2, ..., I_s$ there are $|M| + 2$ segments of length at least $\overline{p} + 2$. Further, since each of the remaining $s - (|M| + 2)$ segments has a length at least 2, we get

$$|C| \geq (\overline{p} + 2)(|M| + 2) + 2(s - |M| - 2)$$
$$= (\overline{p} - 2)|M| + (2\overline{p} + 4|M| + 4) + 2(|A_1| + |A_2| - 2)$$
$$\geq 2|A_1| + 2|A_2| + 4|M| + 2\overline{p}.$$

Observing also that

$$|A_1| + |M| + \overline{p} \geq d(x), \quad |A_2| + |M| + \overline{p} \geq d(y),$$

we have

$$2|A_1| + 2|A_2| + 4|M| + 2\overline{p}$$
$$\geq 2d(x) + 2d(y) - 2\overline{p} \geq 4\delta - 2\overline{p},$$

implying that $|C| \geq 4\delta - 2\overline{p}$.

**Case 2.2**. Either $|A_1| \geq 1, |A_2| = 0$ or $|A_1| = 0, |A_2| \geq 1$.

Assume w.l.o.g. that $|A_1| \geq 1$ and $|A_2| = 0$, i.e. $|N_C(y)| = |M| \geq 2$ and $s = |A_1| + |M|$. It follows that among $I_1, I_2, ..., I_s$ there are $|M| + 1$ segments of length at least $\overline{p} + 2$. Observing also that $|M| + \overline{p} \geq d(y) \geq \delta$, i.e. $2\overline{p} + 4|M| \geq 4\delta - 2\overline{p}$, we get

$$|C| \geq (\overline{p} + 2)(|M| + 1) \geq (\overline{p} - 2)(|M| - 1) + 2\overline{p} + 4|M|$$
$$\geq 2\overline{p} + 4|M| \geq 4\delta - 2\overline{p}. \quad \blacksquare$$

**Proof of Lemma 2**. Let $\xi_1, \xi_2, ..., \xi_s$ be the elements of $N_C(x)$ occuring on $C$ in a consecutive order. Put $I_i = \xi_i \overrightarrow{C} \xi_{i+1}$ $(i = 1, 2, ..., s)$, where $\xi_{s+1} = \xi_1$. To prove $(a1)$, let $L = z\overrightarrow{L}w$ be an intermediate path between elementary segments $I_a$ and $I_b$ with $z \in V(I_a^*)$ and $w \in V(I_b^*)$. Put

$$|\xi_a \overrightarrow{C} z| = d_1, \ |z\overrightarrow{C}\xi_{a+1}| = d_2, \ |\xi_b \overrightarrow{C} w| = d_3, \ |w\overrightarrow{C}\xi_{b+1}| = d_4,$$
$$C' = \xi_a \overrightarrow{P} y \xi_b \overleftarrow{C} z \overrightarrow{L} w \overrightarrow{C} \xi_a.$$

Clearly,
$$|C'| = |C| - d_1 - d_3 + |L| + |P| + 2.$$

Since $C$ is extreme, we have $|C| \geq |C'|$, implying that $d_1 + d_3 \geq \overline{p} + |L| + 2$. By a symmetric argument, $d_2 + d_4 \geq \overline{p} + |L| + 2$. Hence

$$|I_a| + |I_b| = \sum_{i=1}^{4} d_i \geq 2\overline{p} + 2|L| + 4.$$

The proof of $(a1)$ is complete. To proof $(a2)$ and $(a3)$, let $\Upsilon(I_a, I_b) \subseteq E(G)$ and $|\Upsilon(I_a, I_b)| = i$ for some $i \in \{1, 2, 3\}$.

**Case 1.** $i = 1$.

It follows that $\Upsilon(I_a, I_b)$ consists of a unique intermediate edge $L = zw$. By $(a1)$,
$$|I_a| + |I_b| \geq 2\overline{p} + 2|L| + 4 = 2\overline{p} + 6.$$

**Case 2.** $i = 2$.

It follows that $\Upsilon(I_a, I_b)$ consists of two edges $e_1, e_2$. Put $e_1 = z_1 w_1$ and $e_2 = z_2 w_2$, where $\{z_1, z_2\} \subseteq V(I_a^*)$ and $\{w_1, w_2\} \subseteq V(I_b^*)$.

**Case 2.1.** $z_1 \neq z_2$ and $w_1 \neq w_2$.

Assume w.l.o.g. that $z_1$ and $z_2$ occur in this order on $I_a$.

**Case 2.1.1.** $w_2$ and $w_1$ occur in this order on $I_b$.

Put
$$|\xi_a \overrightarrow{C} z_1| = d_1, \ |z_1 \overrightarrow{C} z_2| = d_2, \ |z_2 \overrightarrow{C} \xi_{a+1}| = d_3,$$
$$|\xi_b \overrightarrow{C} w_2| = d_4, \ |w_2 \overrightarrow{C} w_1| = d_5, \ |w_1 \overrightarrow{C} \xi_{b+1}| = d_6,$$
$$C' = \xi_a \overrightarrow{C} z_1 w_1 \overleftarrow{C} w_2 z_2 \overrightarrow{C} \xi_b x \overrightarrow{P} y \xi_{b+1} \overrightarrow{C} \xi_a.$$

Clearly,
$$|C'| = |C| - d_2 - d_4 - d_6 + |\{e_1\}| + |\{e_2\}| + |P| + 2$$
$$= |C| - d_2 - d_4 - d_6 + \overline{p} + 4.$$

Since $C$ is extreme, $|C| \geq |C'|$, implying that $d_2 + d_4 + d_6 \geq \overline{p} + 4$. By a symmetric argument, $d_1 + d_3 + d_5 \geq \overline{p} + 4$. Hence
$$|I_a| + |I_b| = \sum_{i=1}^{6} d_i \geq 2\overline{p} + 8.$$

**Case 2.1.2.** $w_1$ and $w_2$ occur in this order on $I_b$.

Putting
$$C' = \xi_a \overrightarrow{C} z_1 w_1 \overrightarrow{C} w_2 z_2 \overrightarrow{C} \xi_b x \overrightarrow{P} y \xi_{b+1} \overrightarrow{C} \xi_a,$$
we can argue as in Case 2.1.1.

**Case 2.2.** Either $z_1 = z_2$, $w_1 \neq w_2$ or $z_1 \neq z_2$, $w_1 = w_2$.

Assume w.l.o.g. that $z_1 \neq z_2$, $w_1 = w_2$ and $z_1, z_2$ occur in this order on $I_a$.
Put
$$|\xi_a \overrightarrow{C} z_1| = d_1, \ |z_1 \overrightarrow{C} z_2| = d_2, \ |z_2 \overrightarrow{C} \xi_{a+1}| = d_3,$$
$$|\xi_b \overrightarrow{C} w_1| = d_4, \ |w_1 \overrightarrow{C} \xi_{b+1}| = d_5,$$
$$C' = \xi_a x \overrightarrow{P} y \xi_b \overleftarrow{C} z_1 w_1 \overrightarrow{C} \xi_a,$$
$$C'' = \xi_a \overrightarrow{C} z_2 w_1 \overleftarrow{C} \xi_{a+1} x \overrightarrow{P} y \xi_{b+1} \overrightarrow{C} \xi_a.$$

Clearly,
$$|C'| = |C| - d_1 - d_4 + |\{e_1\}| + |P| + 2 = |C| - d_1 - d_4 + \overline{p} + 3,$$
$$|C''| = |C| - d_3 - d_5 + |\{e_2\}| + |P| + 2 = |C| - d_3 - d_5 + \overline{p} + 3.$$
Since $C$ is extreme, $|C| \geq |C'|$ and $|C| \geq |C''|$, implying that
$$d_1 + d_4 \geq \overline{p} + 3, \ d_3 + d_5 \geq \overline{p} + 3.$$
Hence,
$$|I_a| + |I_b| = \sum_{i=1}^{5} d_i \geq d_1 + d_3 + d_4 + d_5 + 1 \geq 2\overline{p} + 7.$$

**Case 3**. $i = 3$.

It follows that $\Upsilon(I_a, I_b)$ consists of three edges $e_1, e_2, e_3$. Let $e_i = z_i w_i$ ($i = 1, 2, 3$), where $\{z_1, z_2, z_3\} \subseteq V(I_a^*)$ and $\{w_1, w_2, w_3\} \subseteq V(I_b^*)$. If there are two independent edges among $e_1, e_2, e_3$ then we can argue as in Case 2.1. Otherwise, we can assume w.l.o.g. that $w_1 = w_2 = w_3$ and $z_1, z_2, z_3$ occur in this order on $I_a$. Put
$$|\xi_a \overrightarrow{C} z_1| = d_1, \ |z_1 \overrightarrow{C} z_2| = d_2, \ |z_2 \overrightarrow{C} z_3| = d_3,$$
$$|z_3 \overrightarrow{C} \xi_{a+1}| = d_4, \ |\xi_b \overrightarrow{C} w_1| = d_5, \ |w_1 \overrightarrow{C} \xi_{b+1}| = d_6,$$
$$C' = \xi_a x \overrightarrow{P} y \xi_b \overleftarrow{C} z_1 w_1 \overrightarrow{C} \xi_a,$$
$$C'' = \xi_a \overrightarrow{C} z_3 w_1 \overleftarrow{C} \xi_{a+1} x \overrightarrow{P} y \xi_{b+1} \overrightarrow{C} \xi_a.$$
Clearly,
$$|C'| = |C| - d_1 - d_5 + |\{e_1\}| + \overline{p} + 2,$$
$$|C''| = |C| - d_4 - d_6 + |\{e_3\}| + \overline{p} + 2.$$
Since $C$ is extreme, we have $|C| \geq |C'|$ and $|C| \geq |C''|$, implying that
$$d_1 + d_5 \geq \overline{p} + 3, \ d_4 + d_6 \geq \overline{p} + 3.$$
Hence,
$$|I_a| + |I_b| = \sum_{i=1}^{6} d_i \geq d_1 + d_4 + d_5 + d_6 + 2 \geq 2\overline{p} + 8. \quad \blacksquare$$

**Proof of Lemma 3**. Choose a longest path $P = x_1 \overrightarrow{P} x_2$ in $G \backslash C$ so as to maximize $|N_C(x_1)|$. Let $y_1, ..., y_t$ be the elements of $N_P^+(x_2)$ occuring on $P$ in a consecutive order. Put
$$P_i = x_1 \overrightarrow{P} y_i^- x_2 \overleftarrow{P} y_i \ (i = 1, ..., t), \quad H = G[V(y_1^- \overrightarrow{P} x_2)].$$
Since $P_i$ is a longest path in $G \backslash C$ for each $i \in \{1, ..., t\}$, we can assume w.l.o.g. that $P$ is chosen such that $|V(H)|$ is maximum. It follows in particular that

$N_P(y_i) \subseteq V(H)$ $(i = 1, ..., t)$.

**Case 1**. $|N_C(x_1)| = 0$.

Since $|N_C(x_1)|$ is maximum, we have $|N_C(y_i)| = 0$ $(i = 1, ..., t)$, implying that $N(y_i) \subseteq V(H)$ and $d_H(y_i) = d(y_i) \geq \delta$ $(i = 1, ..., t)$. Further, since $y_t = x_2$, we have $d_P(x_2) \geq \delta$, that is $t \geq \delta$. By Theorem C, for each distinct $u, v \in V(H)$, there is a path in $H$ of length at least $\delta$, connecting $u$ and $v$. Since $H$ and $C$ are connected by at least $\kappa$ vertex disjoint paths, we have $|C| \geq \kappa(\delta + 2)$.

**Case 2**. $|N_C(x_1)| = 1$.

Since $|N_C(x_1)|$ is maximum, we have $|N_C(y_i)| \leq 1$ $(i = 1, ..., t)$, implying that $|N_H(y_i)| \geq \delta - 1$ $(i = 1, ..., t)$, where $t \geq \delta - 1$. By Theorem C, $|C| \geq \kappa(\delta + 1)$.

**Case 3**. $|N_C(x_1)| \geq 2$.

If $|N_C(y_i)| \geq 2$ for some $i \in \{1, ..., t\}$ then we are done. Otherwise $|N_C(y_i)| \leq 1$ $(i = 1, ..., t)$ and, as in Case 2, $|C| \geq \kappa(\delta + 1)$. ∎

**Proof of Theorem 1**. If $\kappa \leq 2$ then $\tau \leq 1$, contradicting the hypothesis. Let $\kappa \geq 3$. Next, if $c \geq 2\delta + 4$ then we are done. So, we can assume that

$$c \leq 2\delta + 3. \tag{1}$$

Let $C$ be a longest cycle in $G$ and $P = x_1 \overrightarrow{P} x_2$ a longest path in $G \backslash C$ of length $\overline{p}$. If $|V(P)| \leq 0$ then $C$ is a Hamilton cycle and we are done. Let $|V(P)| \geq 1$. Put $X = N_C(x_1) \cup N_C(x_2)$ and let $\xi_1, ..., \xi_s$ be the elements of $X$ occuring on $C$ in a consecutive order. Put

$$I_i = \xi_i \overrightarrow{C} \xi_{i+1}, \ I_i^* = \xi_i^+ \overrightarrow{C} \xi_{i+1}^- \ (i = 1, ..., s),$$

where $\xi_{s+1} = \xi_1$.

**Claim 1**. Let $N_C(x_1) = N_C(x_2)$ and let $\xi_a, \xi_b$ be two distinct elements of $X$. If either $|\xi_a \overrightarrow{C} y| + |\xi_b \overrightarrow{C} z| \leq \overline{p} + 2$ or $|y \overrightarrow{C} \xi_{a+1}| + |z \overrightarrow{C} \xi_{b+1}| \leq \overline{p} + 2$ for some $y \in V(I_a^*)$ and $z \in V(I_b^*)$, then $yz \notin E(G)$.

**Proof**. Assume the contrary, that is $yz \in E(G)$. If $|\xi_a \overrightarrow{C} y| + |\xi_b \overrightarrow{C} z| \leq \overline{p} + 2$ then

$$|\xi_a x_1 \overrightarrow{P} x_2 \xi_b \overleftarrow{C} y z \overrightarrow{C} \xi_a| = |C| - |\xi_a \overrightarrow{C} y| - |\xi_b \overrightarrow{C} z| + \overline{p} + 3 \geq |C| + 1,$$

a contradiction. By a symmetric argument, we reach a contradiction when $|y \overrightarrow{C} \xi_{a+1}| + |z \overrightarrow{C} \xi_{b+1}| \leq \overline{p} + 2$. △

**Claim 2**. Let $N_C(x_1) = N_C(x_2)$ and let $\xi_a, \xi_b, \xi_f$ be distinct elements of $X$, occuring on $\overrightarrow{C}$ in a consecutive order. If $\xi_a^- \xi_b^+ \in E(G)$ and $|\xi_f \overrightarrow{C} y| \leq \overline{p} + 1$ for some $y \in V(I_f^*)$, then $y\xi_a, y\xi_b \notin E(G)$.

**Proof**. Assume the contrary. If $y\xi_a \in E(G)$ then

$$|\xi_f x_1 \overrightarrow{P} x_2 \xi_b \overleftarrow{C} \xi_a y \overrightarrow{C} \xi_a^- \xi_b^+ \overrightarrow{C} \xi_f| = |C| - |\xi_f \overrightarrow{C} y| + \overline{p} + 2 \geq |C| + 1,$$

a contradiction. If $y\xi_b \in E(G)$ then

$$|\xi_f x_1 \overrightarrow{P} x_2 \xi_a \overrightarrow{C} \xi_b y \overrightarrow{C} \xi_a^- \xi_b^+ \overrightarrow{C} \xi_f| \geq |C| + 1,$$

a contradiction. $\Delta$

**Case 1.** $\overline{p} = 0$.

It follows that $P = x_1$ and $s = d(x_1) \geq \delta \geq 3$. Assume first that $s \geq \delta + 1$. If $\Upsilon(I_1, ..., I_s) = \emptyset$ then $G\backslash\{\xi_1, ..., \xi_s\}$ has at least $s+1$ components, contradicting the fact that $\tau > 1$. Otherwise $\Upsilon(I_a, I_b) \neq \emptyset$ for some distinct $a, b \in \{1, ..., s\}$. By Lemma 2, $|I_a| + |I_b| \geq 6$. Since $C$ is extreme, we have $|I_i| \geq 2$ ($i = 1, ..., s$) and therefore, $c \geq 6 + 2(s-2) \geq 2\delta + 4$, contradicting (1). So, $s = \delta$.

The next claim can be derived from (1) and Lemma 2 easily.

**Claim 3.** (1) $|I_i| + |I_j| \leq 7$ for each distinct $i, j \in \{1, ..., s\}$.
(2) If $|I_a| + |I_b| = 7$ for some distinct $a, b \in \{1, ..., s\}$ then $|I_i| = 2$ for each $i \in \{1, ..., s\}\backslash\{a, b\}$.
(3) If $|I_a| = 5$ for some $a \in \{1, ..., s\}$ then $|I_i| = 2$ for each $i \in \{1, ..., s\}\backslash\{a\}$.
(4) There are at most three segments of length at least 3.
(5) If $|I_a| \geq 3$, $|I_b| \geq 3$, $|I_f| \geq 3$ for some distinct $a, b, f \in \{1, ..., s\}$ then $|I_a| = |I_b| = |I_f| = 3$.

If $\Upsilon(I_1, ..., I_s) = \emptyset$ then $G\backslash\{\xi_1, ..., \xi_s\}$ has at least $s+1$ components, contradicting the fact that $\tau > 1$. Otherwise $\Upsilon(I_i, I_j) \neq \emptyset$ for some distinct $i, j \in \{1, ..., s\}$. Choose $a, b \in \{1, ..., s\}$ such that $\Upsilon(I_a, I_b) \neq \emptyset$ and $|I_a| + |I_b|$ is maximum. By definition, there is an intermediate path $L$ between $I_a$ and $I_b$. If $|L| \geq 2$ then by Lemma 2,

$$|I_a| + |I_b| \geq 2\overline{p} + 2|L| + 4 \geq 8,$$

contradicting Claim 3(1). Otherwise $|L| = 1$ and therefore,

$$\Upsilon(I_1, ..., I_s) \subseteq E(G).$$

By Lemma 2, $|I_a| + |I_b| \geq 2\overline{p} + 6 = 6$. Combining this with Claim 3(1), we have

$$6 \leq |I_a| + |I_b| \leq 7.$$

Let $L = yz$, where $y \in V(I_a^*)$ and $z \in V(I_b^*)$.

**Case 1.1.** $|I_a| + |I_b| = 6$.

Since $|I_i| \geq 2$ ($i = 1, ..., s$), we can assume w.l.o.g. that either $|I_a| = 2$, $|I_b| = 4$ or $|I_a| = |I_b| = 3$.

**Case 1.1.1.** $|I_a| = 2$ and $|I_b| = 4$.

Put $I_a = \xi_a w_1 \xi_{a+1}$ and $I_b = \xi_b w_2 w_3 w_4 \xi_{b+1}$. Since $|I_a| + |I_b|$ is extreme, we have $|I_i| = 2$ for each $i \in \{1, ..., s\}\backslash\{b\}$. Clearly, $y = w_1$. By Claim 1, $z = w_3$ and

$\Upsilon(I_a, I_b) = \{w_1 w_3\}$. If $\Upsilon(I_1, ..., I_s) = \{w_1 w_3\}$ then $G \backslash \{\xi_1, ..., \xi_s, w_3\}$ has at least $s+1$ components, contradicting the fact that $\tau > 1$. Otherwise $\Upsilon(I_f, I_g) \neq \emptyset$ for some distinct $f, g \in \{1, ..., s\}$ with $\{f, g\} \neq \{a, b\}$. If $\{f, g\} \cap \{a, b\} = \emptyset$ then by Lemma 2, $|I_f| + |I_g| \geq 6$ and therefore,

$$c = \sum_{i \in \{a,b,f,g\}} |I_i| + \sum_{i \in \{1,2,...,s\} \backslash \{a,b,f,g\}} |I_i| \geq 12 + 2(s-4) = 2\delta + 4,$$

contradicting (1). Let $\{f, g\} \cap \{a, b\} \neq \emptyset$. If $f = a$ then clearly $g \neq b$ and by Lemma 2, $|I_a| + |I_g| \geq 6$, implying that $|I_g| \geq 4$. But then $|I_b| + |I_g| \geq 8$, contradicting Claim 3(1). Now let $f \neq a$ and $g = b$. By Lemma 2, $|I_b| + |I_f| \geq 6$. Since $|I_a| + |I_b|$ is extreme, we have $|I_b| + |I_f| = 6$, which yields $|I_f| = 2$. Put $I_f = \xi_f w_5 \xi_{f+1}$. Let $y_1 z_1$ be an intermediate edge between $I_f$ and $I_b$. By Claim 1, $y_1 = w_5$ and $z_1 = w_3$. Recalling that $|I_i| = 2$ for each $i \in \{1, ..., s\} \backslash \{b\}$, we conclude that $w_3$ belongs to all intermediate edges in $\Upsilon(I_1, ..., I_s)$. Then $G \backslash \{\xi_1, ..., \xi_s, w_3\}$ has at least $s+1$ components, contradicting the fact that $\tau > 1$.

**Case 1.1.2**. $|I_a| = |I_b| = 3$.

Put $I_a = \xi_a w_1 w_2 \xi_{a+1}$ and $I_b = \xi_b w_3 w_4 \xi_{b+1}$. Assume w.l.o.g. that $y = w_2$. By Claim 1, $z = w_3$ and $\Upsilon(I_a, I_b) = \{w_2 w_3\}$. If $\Upsilon(I_1, ..., I_s) = \{w_2 w_3\}$ then $G \backslash \{\xi_1, ..., \xi_s, w_2\}$ has at least $s+1$ components, contradicting the fact that $\tau > 1$. Otherwise $\Upsilon(I_f, I_g) \neq \emptyset$ for some distinct $f, g \in \{1, ..., s\}$ with $\{f, g\} \neq \{a, b\}$. If $\{f, g\} \cap \{a, b\} = \emptyset$ then by Lemma 2, $|I_f| + |I_g| \geq 6$ and, as in Case 1.1.1, $c \geq 12 + 2(s-4) \geq 2\delta + 4$, contradicting (1). Let $\{f, g\} \cap \{a, b\} \neq \emptyset$. Assume w.l.o.g. that $f = a$ and $g \neq b$. By Lemma 2, $|I_a| + |I_g| \geq 6$, that is $|I_g| \geq 3$. By Claim 3(5), $|I_g| = 3$. Put $I_g = \xi_g w_5 w_6 \xi_{g+1}$. Let $y_1 z_1$ be an intermediate edge with $y_1 \in V(I_a^*)$ and $z_1 \in V(I_g^*)$.

**Case 1.1.2.1**. $g \in V(\xi_{b+1}^+ \overrightarrow{C} \xi_a^-)$.

If $y_1 = w_1$ then by Claim 1, $z_1 = w_6$ and

$$\xi_a w_1 w_6 \overleftarrow{C} w_3 w_2 \overrightarrow{C} \xi_b x_1 \xi_{g+1} \overrightarrow{C} \xi_a$$

is longer than $C$, a contradiction. Let $y_1 = w_2$. By Claim 1, $z_1 = w_5$ and therefore, $\Upsilon(I_a, I_g) = \{w_2 w_5\}$.

**Case 1.1.2.1.1**. $N(w_1) \subseteq V(C)$.

By Claim 2, $w_1 \xi_b \notin E(G)$ and $w_1 \xi_g \notin E(G)$. Further, if

$$N(w_1) \subseteq \{\xi_1, ..., \xi_s, w_2\} \backslash \{\xi_b, \xi_g\}$$

then $|N(w_1)| \leq s - 1 = \delta - 1$, a contradiction. Otherwise, $w_1 z_2 \in E(G)$ for some $z_2 \in V(I_h^*)$, where $h \notin \{a, b, g\}$. By Lemma 2, $|I_a| + |I_h| \geq 6$, implying that $|I_h| \geq 3$, which contradicts Claim 3(4).

**Case 1.1.2.1.2**. $N(w_1) \not\subseteq V(C)$.

It follows that $w_1x_2 \in E(G)$ for some $x_2 \in V(G\backslash C)$. Since $\overline{p} = 0$ and $C$ is extreme, $x_2 \neq x_1$ and $N(x_2) \subseteq V(C)$. By the same reason, $x_2\xi_a \notin E(G)$ and $x_2w_2 \notin E(G)$. By Claim 2, $x_2\xi_b \notin E(G)$. If

$$N(x_2) \subseteq \{\xi_1, ..., \xi_s, w_1\}\backslash\{\xi_a, \xi_b\}$$

then $|N(x_2)| \leq s - 1 = \delta - 1$, a contradiction. Otherwise $x_2z_2 \in E(G)$ for some $z_2 \in V(I_h^*)$, where $h \neq a$. But then $I_a^*$ and $I_h^*$ are connected by $w_1x_2z_2$, contradicting the fact that $\Upsilon(I_1, ..., I_s) \subseteq E(G)$.

**Case 1.1.2.2.** $g \in V(\xi_{a+1}^+ \overrightarrow{C} \xi_b^-)$.

If $y_1 = w_2$ then by Claim 1, $z_1 = w_5$ and we can argue as in Case 1.1.2.1. Let $z_1 = w_1$. By Claim 1, $z_2 = w_6$ and $w_4w_6 \notin E(G)$. Further, by Claim 2, $w_4\xi_{a+1} \notin E(G)$ and $w_4\xi_b \notin E(G)$. Using Claim 3(4), we have $|I_i| = 2$ for each $i \in \{1, ..., s\}\backslash\{a, b, g\}$. By Lemma 2, $N(w_4) \cap V(I_i^*) = \emptyset$ for each $i \in \{1, ..., s\}\backslash\{a, b, g\}$.

**Case 1.1.2.2.1.** $N(w_4) \subseteq V(C)$.

It follows that

$$N(w_4) \subseteq \{\xi_1, ..., \xi_s, w_3, w_5\}\backslash\{\xi_{a+1}, \xi_b\}.$$

Since $|N(w_4)| \geq \delta = s$, we have $w_4w_5 \in E(G)$.

**Case 1.1.2.2.1.1.** $s \geq 4$.

Since $|I_i| = 2$ for each $i \in \{1, ..., s\}\backslash\{a, b, g\}$, we can assume w.l.o.g. that $|I_{a-1}| = 2$. Put $I_{a-1} = \xi_{a-1}w_7\xi_a$. Assume first that $N(w_7) \not\subseteq V(C)$, that is $w_7x_2 \in E(G)$ for some $x_2 \in V(G\backslash C)$. Since $C$ is extreme and $\Upsilon(I_1, ..., I_s) \subseteq E(G)$, we have $x_2 \neq x_1$ and

$$N(x_2) \subseteq \{\xi_1, ..., \xi_s, w_7\}\backslash\{\xi_{a-1}, \xi_a, \},$$

contradicting the fact that $|N(x_2)| \geq \delta = s$. Now assume that $N(w_7) \subseteq V(C)$. By Claim 2, $w_7\xi_{a+1} \notin E(G)$. Since $|I_{a-1}| = 2$ and $|I_i| \leq 3$ for each $i \in \{1, ..., s\}$, we have by Lemma 2, $N(w_7) \cap V(I_i^*) = \emptyset$ for each $i \in \{1, ..., s\}\backslash\{a - 1\}$. So, $N(w_7) \subseteq \{\xi_1, ..., \xi_s\}\backslash\{\xi_{a+1}\}$, contradicting the fact that $|N(w_7)| \geq \delta = s$.

**Case 1.1.2.2.1.2.** $s = 3$.

Put $C = \xi_1w_1w_2\xi_2w_3w_4\xi_3w_5w_6\xi_1$. Assume first that $N(w_i) \not\subseteq V(C)$ for some $i \in \{1, 2, ..., 6\}$, say $i = 1$. This means that $w_1x_2 \in E(G)$ for some $x_2 \in V(G\backslash C)$. Since $C$ is extreme, $x_2 \neq x_1$ and $x_2\xi_1, x_2w_2 \notin E(G)$. Further, since $\Upsilon(I_1, I_2, I_3) \subseteq E(G)$, we have $N(x_2) \subseteq \{\xi_2, \xi_3, w_1\}$. On the other hand, since $|N(x_2)| \geq \delta \geq 3$, we have $N(x_2) = \{\xi_2, \xi_3, w_1\}$. By Claim 2, $x_2\xi_2 \notin E(G)$, a contradiction. Now assume that $N(w_i) \subseteq V(C)$ ($i = 1, ..., 6$). If $V(G\backslash C) \neq \{x_1\}$, then choose $x_2 \in V(G\backslash C)$ such that $x_2 \neq x_1$. Since $N(w_i) \subseteq V(C)$ ($i = 1, ..., 6$), we have $N(x_2) = N(x_1)$. But then $G\backslash\{\xi_1, \xi_2, \xi_3\}$ has at least three components, contradicting the fact that $\tau > 1$. Finally, if $V(G\backslash C) = \{x_1\}$, then $G$ is the Petersen graph.

**Case 1.1.2.2.2**. $N(w_4) \not\subseteq V(C)$.

It follows that $w_4 x_2 \in E(G)$ for some $x_2 \in V(G\backslash C)$. Since $C$ is extreme and $\Upsilon(I_1, ..., I_s) \subseteq E(G)$, we have $x_2 \neq x_1$, $x_2 \xi_{b+1} \notin E(G)$ and $N(x_2) \cap V(I_i^*) = \emptyset$ for each $i \in \{1, ..., s\}\backslash\{b\}$. So, $N(x_2) \subseteq \{\xi_1, ..., \xi_s, w_4\}\backslash\{\xi_{b+1}\}$, implying that $x_2 \xi_b \in E(G)$ which contradicts Claim 2.

**Case 1.2**. $|I_a| + |I_b| = 7$.

By Claim 3(2), $|I_i| = 2$ for each $i \in \{1, ..., s\}\backslash\{a, b\}$. By the hypothesis, either $|I_a| = 2$, $|I_b| = 5$ or $|I_a| = 3$, $|I_b| = 4$.

**Case 1.2.1**. $|I_a| = 2$, $|I_b| = 5$.

Put $I_a = \xi_a w_1 \xi_{a+1}$ and $I_b = \xi_b w_2 w_3 w_4 w_5 \xi_{b+1}$. Clearly, $y = w_1$. By Claim 1, $z \in \{w_3, w_4\}$. Further, if $\{w_1 w_3, w_1 w_4\} \subseteq E(G)$ then

$$\xi_a x_1 \xi_{a+1} \overrightarrow{C} w_3 w_1 w_4 \overrightarrow{C} \xi_a$$

is longer than $C$, a contradiction. Therefore, we can assume w.l.o.g. that $w_1 w_3 \in E(G)$ and $\Upsilon(I_a, I_b) = \{w_1 w_3\}$. By Claim 3(3), $|I_i| = 2$ for each $i \in \{1, ..., s\}\backslash\{b\}$. By Lemma 2, each intermediate edge has one end in $V(I_b^*)$. If $\Upsilon(I_1, ..., I_s) = \{w_1 w_3\}$ then $G\backslash\{\xi_1, ..., \xi_s, w_3\}$ has at least $s + 1$ components, contradicting the fact that $\tau > 1$. Otherwise $\Upsilon(I_b, I_g) \neq \emptyset$ for some $g \in \{1, ..., s\}\backslash\{a, b\}$. Since $|I_g| = 2$, we can set $I_g = \xi_g w_6 \xi_{g+1}$. As above, either $w_6 w_3 \in E(G)$, $w_6 w_4 \notin E(G)$ or $w_6 w_3 \notin E(G)$, $w_6 w_4 \in E(G)$. Assume that $w_6 w_4 \in E(G)$. If $\xi_g \in V(\xi_{b+1}^+ \overrightarrow{C} \xi_a^-)$ then

$$\xi_a w_1 w_3 \overleftarrow{C} \xi_{a+1} x_1 \xi_g \overleftarrow{C} w_4 w_6 \overrightarrow{C} \xi_a$$

is longer than $C$, a contradiction. If $\xi_g \in V(\xi_{a+1}^+ \overrightarrow{C} \xi_b^-)$ then

$$\xi_a x_1 \xi_{g+1} \overrightarrow{C} w_3 w_1 \overrightarrow{C} w_6 w_4 \overrightarrow{C} \xi_a$$

is longer than $C$, a contradiction. Now assume that $w_6 w_4 \notin E(G)$, implying that $w_6 w_3 \in E(G)$. This means that $w_3$ belongs to each intermediate edge in $\Upsilon(I_1, ..., I_s)$. But then $G\backslash\{\xi_1, ..., \xi_s, w_3\}$ has at least $s + 1$ components, contradicting the fact that $\tau > 1$.

**Case 1.2.2**. $|I_a| = 3$, $|I_b| = 4$.

Put $I_a = \xi_a w_1 w_2 \xi_{a+1}$ and $I_b = \xi_b w_3 w_4 w_5 \xi_{b+1}$. Assume w.l.o.g. that $y = w_2$. By Claim 1, $z \in \{w_3, w_4\}$.

**Case 1.2.2.1**. $w_2 w_3 \in E(G)$.

Assume first that $N(w_1) \not\subseteq V(C)$, that is $w_1 x_2 \in E(G)$ for some $x_2 \in V(G\backslash C)$. Since $C$ is extreme, $x_2 \neq x_1$ and $x_2 \xi_a \notin E(G)$, $x_2 w_2 \notin E(G)$. By Claim 2, $x_2 \xi_{a+1} \notin E(G)$. Recalling that $C$ is extreme and $\Upsilon(I_1, ..., I_s) \subseteq E(G)$, we have

$$N(x_2) \subseteq \{\xi_1, ..., \xi_s, w_1\}\backslash\{\xi_a, \xi_{a+1}\},$$

contradicting the fact that $|N(x_2)| \geq \delta = s$. Now assume that $N(w_1) \subseteq V(C)$. By Claim 1, $w_1\xi_{a+1} \notin E(G)$, $w_1\xi_b \notin E(G)$ and $w_1w_3 \notin E(G)$. Further, if $N(w_1) \cap \{w_4, w_5\} \neq \emptyset$ then there are two independent intermediate edges between $I_a$ and $I_b$. By Lemma 2, $|I_a| + |I_b| \geq 8$, contradicting Claim 3(1). Hence, $N(w_1) \cap \{w_4, w_5\} = \emptyset$. Finally, since $|I_i| = 2$ for each $i \in \{1, ..., s\}\backslash\{a, b\}$, we have $N(w_1) \cap V(I_i^*) = \emptyset$ for each $i \in \{1, ..., s\}\backslash\{a\}$. So,

$$N(w_1) \subseteq \{\xi_1, ..., \xi_s, w_2\}\backslash\{\xi_{a+1}, \xi_b\},$$

contradicting the fact that $|N(w_1)| \geq \delta = s$ when $\xi_{a+1} \neq \xi_b$. Let $\xi_{a+1} = \xi_b$. Assume w.l.o.g. that $a = 1$ and $b = 2$. If $s = 2$ then clearly $\tau \leq 1$, contradicting the hypothesis. Let $s \geq 3$. Recalling that $|I_i| = 2$ for each $i \in \{3, ..., s\}$, we can set $I_3 = \xi_3 w_7 \xi_4$. If $N(w_7) \not\subseteq V(C)$, that is $w_7 x_2 \in E(G)$ for some $x_2 \in V(G\backslash C)$, then $x_2 \neq x_1$ and

$$N(x_2) \subseteq \{\xi_1, ..., \xi_s, w_7\}\backslash\{\xi_3, \xi_4\},$$

contradicting the fact that $|N(x_2)| \geq \delta = s$. Let $N(w_7) \subseteq V(C)$. By Claim 2, $w_7\xi_2 \notin E(G)$. Hence, $N(w_7) \subseteq \{\xi_1, ..., \xi_s\}\backslash\{\xi_2\}$, contradicting the fact that $|N(w_7)| \geq s$.

**Case 1.2.2.2.** $w_2w_4 \in E(G)$.

If $w_2w_3 \in E(G)$ then we can argue as in Case 1.2.2.1. Hence, we can assume that $\Upsilon(I_a, I_b) = \{w_2w_4\}$. If $\Upsilon(I_1, ..., I_s) = \{w_2w_4\}$ then clearly $\tau \leq 1$, contradicting the hypothesis. Let $\Upsilon(I_1, ..., I_s) \neq \{w_2w_4\}$. Since $|I_a| = 3$ and $|I_i| = 2$ for each $i \in \{1, ..., s\}\backslash\{a, b\}$, we can state by Lemma 2 that each intermediate edge has one end in $V(I_b^*)$. Let $y_1z_1 \in E(G)$ for some $y_1 \in V(I_g^*)$ and $z_1 \in V(I_b^*)$, where $g \in \{1, ..., s\}\backslash\{a, b\}$. Since $|I_g| = 2$, we can set $I_g = \xi_g w_6 \xi_{g+1}$. Clearly $y_1 = w_6$. By Claim 1, $z_1 = w_4$. This means that $w_4$ belongs to all intermediate edges. Then clearly $\tau \leq 1$, contradicting the hypothesis.

**Case 2.** $\overline{p} = 1$.
Since $\delta \geq \kappa \geq 3$, we have $|N_C(x_i)| \geq \delta - \overline{p} = \delta - 1 \geq 2$ $(i = 1, 2)$.

**Case 2.1.** $N_C(x_1) \neq N_C(x_2)$.
It follows that $\max\{\sigma_1, \sigma_2\} \geq 1$, where

$$\sigma_1 = |N_C(x_1)\backslash N_C(x_2)|, \quad \sigma_2 = |N_C(x_2)\backslash N_C(x_1)|.$$

If $\max\{\sigma_1, \sigma_2\} \geq 2$ then by Lemma 1, $c \geq 3\delta + 1 \geq 2\delta + 4$, contradicting (1). Let $\max\{\sigma_1, \sigma_2\} = 1$. This implies $s \geq \delta$ and $|I_i| \geq 3$ $(i = 1, ..., s)$. If $s \geq \delta + 1$ then $c \geq 3s \geq 3\delta + 3 > 2\delta + 4$, again contradicting (1). Let $s = \delta$, that is $|I_i| = 3$ $(i = 1, ..., s)$. By Lemma 2, $\Upsilon(I_1, ..., I_s) = \emptyset$, contradicting the fact that $\tau > 1$.

**Case 2.2.** $N_C(x_1) = N_C(x_2)$.
Clearly, $s = |N_C(x_1)| \geq \delta - \overline{p} = \delta - 1$. If $s \geq \delta$ then $c \geq 3s \geq 3\delta$ and we can argue as in Case 2.1. Let $s = \delta - 1$.
The following can be derived from (1) and Lemma 2 easily.

**Claim 4**. (1) $|I_i| + |I_j| \leq 9$ for each distinct $i, j \in \{1, ..., s\}$.

(2) If $|I_a| + |I_b| = 9$ for some distinct $a, b \in \{1, ..., s\}$ then $|I_i| = 3$ for each $i \in \{1, ..., s\}\backslash\{a, b\}$.

(3) If $|I_a| = 6$ for some $a \in \{1, ..., s\}$ then $|I_i| = 3$ for each $i \in \{1, ..., s\}\backslash\{a\}$.

(4) There are at most three segments of length at least 4.

(5) If $|I_a| \geq 4$, $|I_b| \geq 4$, $|I_f| \geq 4$ for some distinct $a, b, f \in \{1, 2, ..., s\}$ then $|I_a| = |I_b| = |I_f| = 4$.

If $\Upsilon(I_1, ..., I_s) = \emptyset$ then clearly, $\tau \leq 1$, contradicting the hypothesis. Otherwise $\Upsilon(I_a, I_b) \neq \emptyset$ for some distinct $a, b \in \{1, ..., s\}$. By definition, there is an intermediate path $L$ between $I_a$ and $I_b$. If $|L| \geq 2$ then by Lemma 2,

$$|I_a| + |I_b| \geq 2\overline{p} + 2|L| + 4 \geq 10,$$

contradicting Claim 4(1). Otherwise $|L| = 1$ and therefore,

$$\Upsilon(I_1, ..., I_s) \subseteq E(G).$$

By Lemma 2, $|I_a| + |I_b| \geq 2\overline{p} + 6 = 8$. Combining this with Claim 4(1), we have

$$8 \leq |I_a| + |I_b| \leq 9.$$

Let $L = yz$, where $y \in V(I_a^*)$ and $z \in V(I_b^*)$.

**Case 2.2.1**. $|I_a| + |I_b| = 8$.

Since $|I_i| \geq 3$ $(i = 1, ..., s)$, we can assume w.l.o.g. that either $|I_a| = 3$, $|I_b| = 5$ or $|I_a| = |I_b| = 4$.

**Case 2.2.1.1**. $|I_a| = 3$ and $|I_b| = 5$.

Put $I_a = \xi_a w_1 w_2 \xi_{a+1}$ and $I_b = \xi_b w_3 w_4 w_5 w_6 \xi_{b+1}$. Assume w.l.o.g. that $y = w_2$. By Claim 1, $z = w_4$. By the same reason, $N(w_1) \cap V(I_b^*) \subseteq \{w_5\}$. If $w_1 w_5 \in E(G)$ then there exist two independent intermediate edges between $I_a$ and $I_b$, which by Lemma 2 yields $|I_a| + |I_b| \geq 2\overline{p} + 8 = 10$, contradicting Claim 4(1). So, $N(w_1) \cap V(I_b^*) = \emptyset$. Further, if $\Upsilon(I_a, I_f) \neq \emptyset$ for some $f \in \{1, ..., s\}\backslash\{a, b\}$ then by Lemma 2, $|I_a| + |I_f| \geq 2\overline{p} + 6 = 8$, implying that $|I_f| \geq 5$. But then $|I_b| + |I_f| \geq 10$, contradicting Claim 4(1). Hence $\Upsilon(I_a, I_i) = \emptyset$ for each $i \in \{1, ..., s\}\backslash\{a, b\}$. By Claim 2, $w_1 \xi_{a+1} \not\in E(G)$. Thus, if $N(w_1) \subseteq V(C)$ then

$$N(w_1) \subseteq \{\xi_1, ..., \xi_s, w_2\}\backslash\{\xi_{a+1}\},$$

contradicting the fact that $|N(w_1)| \geq \delta = s + 1$. Now let $N(w_1) \not\subseteq V(C)$ and let $Q = w_1 \overrightarrow{Q} x_3$ be a longest path having only $w_1$ in common with $C$. Clearly, $1 \leq |Q| \leq 2$ and $V(Q) \cap V(P) = \emptyset$. By Claim 2, $x_3 \xi_{a+1} \not\in E(G)$. Further, since $\Upsilon(I_1, ..., I_s) \subseteq E(G)$, we have $N(x_3) \cap V(I_i^*) = \emptyset$ for each $i \in \{1, ..., s\}\backslash\{a\}$. If $|Q| = 1$ then

$$N(x_3) \subseteq \{\xi_1, ..., \xi_s, w_1\}\backslash\{\xi_a, \xi_{a+1}\},$$

contradicting the fact that $|N(x_3)| \geq \delta = s + 1$. If $|Q| = 2$ then

$$N(x_3) \subseteq \{\xi_1, ..., \xi_s, x_3^-, w_1\}\backslash\{\xi_a, \xi_{a+1}\},$$

contradicting the fact that $|N(x_3)| \geq \delta = s+1$.

**Case 2.2.1.2.** $|I_a| = |I_b| = 4$.
Put $I_a = \xi_a w_1 w_2 w_3 \xi_{a+1}$ and $I_b = \xi_b w_4 w_5 w_6 \xi_{b+1}$.

**Case 2.2.1.2.1.** $y \in \{w_1, w_3\}$.
Assume w.l.o.g. that $y = w_3$. By Claim 1, $z = w_4$.

**Claim 5**. $N(w_1) \cup N(w_2) \subseteq V(C)$.
**Proof**. Assume the contrary and let $Q = w_1 \overrightarrow{Q} x_3$ be a longest path having only $w_1$ in common with $C$. Clearly, $1 \leq |Q| \leq 2$ and $V(Q) \cap V(P) = \emptyset$. By Claim 2, $x_3 \xi_{a+1} \notin E(G)$ and $x_3 \xi_b \notin E(G)$. Since $\Upsilon(I_1, ..., I_s) \subseteq E(G)$, we have $N(x_3) \cap V(I_i^*) = \emptyset$ for each $i \in \{1, ..., s\} \backslash \{a\}$. If $|Q| = 1$ then

$$N(x_3) \subseteq \{\xi_1, ..., \xi_s, w_1, w_3\} \backslash \{\xi_a, \xi_{a+1}\},$$

contradicting the fact that $|N(x_3)| \geq \delta = s+1$. If $|Q| = 2$ then

$$N(x_3) \subseteq \{\xi_1, ..., \xi_s, x_3^-, w_1\} \backslash \{\xi_a, \xi_{a+1}\},$$

a contradiction. Similarly, we can reach a contradiction when $N(w_2) \not\subseteq V(C)$. Claim 5 is proved. $\triangle$

**Case 2.2.1.2.1.1.** $\xi_{a+1} \neq \xi_b$.
By Claim 2, $w_1 \xi_{a+1} \notin E(G)$ and $w_1 \xi_b \notin E(G)$. By Claim 1, $w_1 w_4 \notin E(G)$. Moreover, if $N(w_1) \cap V(I_b^*) \neq \emptyset$ then there exist two independent intermediate edges between $I_a$ and $I_b$ which by Lemma 2 yields $|I_a| + |I_b| \geq 2\overline{p} + 8 \geq 10$, contradicting Claim 4(1). Furthermore, if $N(w_1) \cap V(I_i^*) = \emptyset$ for each $i \in \{1, ..., s\} \backslash \{a, b\}$ then by Claim 5,

$$N(w_1) \subseteq \{\xi_1, ..., \xi_s, w_2, w_3\} \backslash \{\xi_{a+1}, \xi_b\},$$

implying that $|N(w_1)| \leq s = \delta - 1$, a contradiction. Otherwise, $w_1 v \in E(G)$, where $v \in V(I_f^*)$ for some $f \in \{1, ..., s\} \backslash \{a, b\}$. By a similar way, it can be shown that $w_2 u \in E(G)$, where $u \in V(I_g^*)$ for some $g \in \{1, ..., s\} \backslash \{a, b\}$. By Lemma 2, $|I_a| + |I_f| \geq 2\overline{p} + 6 = 8$, that is $|I_f| \geq 4$. By Claim 4(5), $|I_f| = 4$. By a symmetric argument, $|I_d| = 4$. Put $I_f = \xi_f w_7 w_8 w_9 \xi_{f+1}$. By Claim 1, $v = w_9$, i.e. $w_1 w_9 \in E(G)$. If $d = f$ then $|\Upsilon(I_a, I_f)| = 2$ and by Lemma 2, $|I_a| + |I_f| \geq 2\overline{p} + 7 = 9$, a contradiction. Otherwise, there are at least four elementary segments of length at least 4, contradicting Claim 4(4).

**Case 2.2.1.2.1.2.** $\xi_{a+1} = \xi_b$.
Assume w.l.o.g. that $a = 1$ and $b = 2$. If $\Upsilon(I_1, I_2, ..., I_s) = \Upsilon(I_1, I_2) = \{w_3 w_4\}$ then clearly, $\tau \leq 1$, a contradiction. Otherwise, there is an intermediate edge $uv$ such that $u \in V(I_1^*) \cup V(I_2^*)$ and $v \in V(I_f^*)$ for some $f \in \{1, 2, ..., s\} \backslash \{1, 2\}$. Assume w.l.o.g. that $u \in V(I_1^*)$. If $u = w_3$ then as above, $\xi_2 = \xi_f$, a contradiction. Let $u \neq w_3$. By Lemma 2, $|I_1| + |I_f| \geq 8$, i.e. $|I_f| \geq 4$.

By Claim 4(5), $|I_f| = 4$. Put $I_f = \xi_c w_7 w_8 w_9 \xi_{f+1}$. If $u = w_1$ then by Claim 1, $v = w_9$ and
$$|\xi_1 w_1 w_9 \overleftarrow{C} w_4 w_3 \xi_2 x_2 x_1 \xi_{f+1} \overrightarrow{C} \xi_1| \geq |C| + 1,$$
a contradiction. If $u = w_2$ then by Claim 1, $v = w_8$ and
$$|\xi_1 w_1 w_2 w_8 \overleftarrow{C} w_4 w_3 \xi_2 x_2 x_1 \xi_{f+1} \overrightarrow{C} \xi_1| \geq |C| + 1,$$
again a contradiction.

**Case 2.2.1.2.2**. $y = w_2$.

By Claim 1, $z = w_5$ and $\Upsilon(I_a, I_b) = \{w_2 w_5\}$. If $|I_i| = 3$ for each $i \in \{1, 2, ..., s\} \backslash \{a, b\}$ then by Lemma 2, $\Upsilon(I_1, I_2, ..., I_s) = \{w_2 w_5\}$ and $\tau \leq 1$, contradicting the hypothesis. Otherwise, $|I_f| \geq 4$ for some $f \in \{1, 2, ..., s\} \backslash \{a, b\}$ and $|I_i| = 3$ for each $i \in \{1, 2, ..., s\} \backslash \{a, b, f\}$. By Claim 4(5), $|I_f| = 4$. Put $I_f = \xi_f w_7 w_8 w_9 \xi_{f+1}$. Clearly, $\Upsilon(I_1, I_2, ..., I_s) = \Upsilon(I_a, I_b, I_f)$. If $\Upsilon(I_a, I_f) = \Upsilon(I_b, I_f) = \emptyset$ then again $\tau \leq 1$, a contradiction. Let $uv \in E(G)$, where $u \in I_a^* \cup I_b^*$ and $v \in V(I_f^*)$. Assume w.l.o.g. that $u \in V(I_a^*)$. If $u \in \{w_1, w_3\}$ then we can argue as in Case 1.2.1.2.1. Let $u = w_2$. By Claim 1, $v = w_8$. If $w_1 w_3 \in E(G)$ then
$$\xi_a x_1 x_2 \xi_b \overleftarrow{C} w_3 w_1 w_2 w_5 \overrightarrow{C} \xi_a$$
is longer than $C$, a contradiction. Let $w_1 w_3 \notin E(G)$. Analogously, $w_4 w_6 \notin E(G)$ and $w_7 w_9 \notin E(G)$. But Then $\{w_1, w_3, w_4, w_6, w_7, w_9\}$ is an independent set of vertices and $G \backslash \{\xi_1, ..., \xi_s, w_2, w_5, w_8\}$ has at least $s + 4$ connected components. Hence $\tau < 1$, contradicting the hypothesis.

**Case 2.2.2**. $|I_a| + |I_b| = 9$.

Since $|I_i| \geq 3$ ($i = 1, ..., s$), we can assume w.l.o.g. that either $|I_a| = 3$, $|I_b| = 6$ or $|I_a| = 4$, $|I_b| = 5$.

**Case 2.2.2.1**. $|I_a| = 3$ and $|I_b| = 6$.

By Claim 4(3), $|I_i| = 3$ for each $i \in \{1, ..., s\} \backslash \{b\}$. Put
$$I_a = \xi_a w_1 w_2 \xi_{a+1}, \quad I_b = \xi_b w_3 w_4 w_5 w_6 w_7 \xi_{b+1}.$$
Since $|I_a| = 3$, we can assume w.l.o.g. that $y = w_2$. By Claim 1, $z \in \{w_4, w_5\}$.

**Case 2.2.2.1.1**. $z = w_4$.

By Claim 1, $w_1 w_4 \notin E(G)$. Next, if $N(w_1) \cap V(I_b^*) \neq \emptyset$ then there are two independent intermediate edges between $I_a$ and $I_b$ and by Lemma 2, $|I_a| + |I_b| \geq 2\overline{p} + 8 = 10$, contradicting Claim 4(1). By Claim 2, $w_1 \xi_{a+1} \notin E(G)$. Finally, by Lemma 2 and Claim 4(3), $N(w_1) \cap V(I_i^*) = \emptyset$ for each $i \in \{1, ..., s\} \backslash \{a, b\}$. So, if $N(w_1) \subseteq V(C)$ then
$$N(w_1) \subseteq \{\xi_1, ..., \xi_s, w_2\} \backslash \{\xi_{a+1}\},$$
contradicting the fact that $|N(w_1)| \geq \delta = s + 1$. Now assume that $N(w_1) \nsubseteq V(C)$. Choose a longest path $Q = w_1 \overrightarrow{Q} x_3$ having only $w_1$ in common with $C$.

Clearly, $V(Q) \cap V(P) = \emptyset$. Since $C$ is extreme, $x_3\xi_a \notin E(G)$ and $x_3x_2 \notin E(G)$. If $x_3\xi_{a+1} \in E(G)$ then

$$\xi_a x_1 x_2 \xi_b \overleftarrow{C} \xi_{a+1} x_3 \overleftarrow{Q} w_1 w_2 w_4 \overrightarrow{C} \xi_a$$

is longer than $C$, a contradiction. Let $x_3\xi_{a+1} \notin E(G)$. If $|Q| = 1$ then

$$N(x_3) \subseteq \{\xi_1, ..., \xi_s, w_1\} \setminus \{\xi_a, \xi_{a+1}\},$$

contradicting the fact that $|N(x_3)| \geq \delta = s+1$. If $|Q| = 2$ then

$$N(x_3) \subseteq \{\xi_1, ..., \xi_s, x_3^-, w_1\} \setminus \{\xi_a, \xi_{a+1}\},$$

contradicting the fact that $|N(x_3)| \geq \delta = s+1$.

**Case 2.2.2.1.2**. $z = w_5$.

If $w_2w_4 \in E(G)$ then we can argue as in Case 2.2.2.1.1. Let $w_2w_4 \notin E(G)$. It means that $w_5$ belongs to all intermediate edges. This implies $\tau \leq 1$, contradicting the hypothesis.

**Case 2.2.2.2**. $|I_a| = 4$ and $|I_b| = 5$.

By Claim 4(2), $|I_i| = 3$ and $\Upsilon(I_a, I_i) = \emptyset$ for each $i \in \{1, ..., s\} \setminus \{a, b\}$. If $\Upsilon(I_b, I_f) \neq \emptyset$ for some $f \in \{1, ..., s\} \setminus \{a, b\}$ then we can argue as in Case 2.2.1.1. Otherwise $\Upsilon(I_1, ..., I_s) = \Upsilon(I_a, I_b)$. If there are two independent edges in $\Upsilon(I_a, I_b)$ then by Lemma 2, $|I_a| + |I_b| \geq 10$, contradicting Claim 4(1). Otherwise $\tau \leq 1$, a contradiction.

**Case 3**. $2 \leq \overline{p} \leq \delta - 3$.

It follows that $|N_C(x_i)| \geq \delta - \overline{p} \geq 3$ $(i = 1, 2)$. If $N_C(x_1) \neq N_C(x_2)$ then by Lemma 1, $|C| \geq 4\delta - 2\overline{p} \geq 3\delta - \overline{p} + 3 \geq 2\delta + 4$, contradicting (1). Hence $N_C(x_1) = N_C(x_2)$, implying that $|I_i| \geq \overline{p} + 2$ $(i = 1, 2, ..., s)$. Clearly, $s \geq |N_C(x_1)| - (|V(P)| - 1) \geq \delta - \overline{p} \geq 3$. If $s \geq \delta - \overline{p} + 1$ then

$$|C| \geq s(\overline{p} + 2) \geq (\delta - \overline{p} + 1)(\overline{p} + 2)$$

$$= (\delta - \overline{p} - 1)(\overline{p} - 1) + 3\delta - \overline{p} + 1 \geq 3\delta - \overline{p} + 3 \geq 2\delta + 4,$$

again contradicting (1). Hence $s = \delta - \overline{p}$. It means that $x_1x_2 \in E(G)$, that is $G[V(P)]$ is hamiltonian. By symmetric arguments, $N_C(y) = N_C(x_1)$ for each $y \in V(P)$. If $\Upsilon(I_1, I_2, ..., I_s) = \emptyset$ then $\tau \leq 1$, contradicting the hypothesis. Otherwise $\Upsilon(I_a, I_b) \neq \emptyset$ for some elementary segments $I_a$ and $I_b$. By definition, there is an intermediate path $L$ between $I_a$ and $I_b$. If $|L| \geq 2$ then by lemma 2,

$$|I_a| + |I_b| \geq 2\overline{p} + 2|L| + 4 \geq 2\overline{p} + 8.$$

Hence

$$|C| = |I_a| + |I_b| + \sum_{i \in \{1,...,s\} \setminus \{a,b\}} |I_i| \geq 2\overline{p} + 8 + (s-2)(\overline{p} + 2)$$

$$= (\delta - \overline{p} - 2)(\overline{p} - 1) + 3\delta - \overline{p} + 2 \geq 3\delta - \overline{p} + 3 \geq 2\delta + 4,$$

contradicting (1). Thus, $|L| = 1$, i.e. $\Upsilon(I_1, I_2, ..., I_s) \subseteq E(G)$. By Lemma 2,

$$|I_a| + |I_b| \geq 2\overline{p} + 2|L| + 4 = 2\overline{p} + 6,$$

which yields

$$|C| = |I_a| + |I_b| + \sum_{i \in \{1,...,s\} \setminus \{a,b\}} |I_i| \geq 2\overline{p} + 6 + (s-2)(\overline{p} + 2)$$

$$= (s-2)(\overline{p} - 2) + (\delta - \overline{p} - 3) + (3\delta - \overline{p} + 1) \geq 3\delta - \overline{p} + 1 \geq 2\delta + 4,$$

contradicting (1).

**Case 4.** $2 \leq \overline{p} = \delta - 2$.

It follows that $|N_C(x_i)| \geq \delta - \overline{p} = 2$ $(i = 1, 2)$. If $N_C(x_1) \neq N_C(x_2)$ then by Lemma 1, $|C| \geq 4\delta - 2\overline{p} = 3\delta - \overline{p} + 2 = 2\delta + 4$, contradicting (1). Hence, $N_C(x_1) = N_C(x_2)$. Clearly, $s = |N_C(x_1)| \geq 2$. Further, if $s \geq 3$ then

$$|C| \geq s(\overline{p} + 2) \geq 3\delta \geq 3\delta - \overline{p} + 2 = 2\delta + 4,$$

again contradicting (1). Hence, $s = 2$. It follows that $x_1 x_2 \in E(G)$, that is $G[V(P)]$ is hamiltonian. By symmetric arguments, $N_C(v) = N_C(x_1) = \{\xi_1, \xi_2\}$ for each $v \in V(P)$. If $\Upsilon(I_1, I_2) = \emptyset$ then clearly, $\tau \leq 1$, contradicting the hypothesis. Otherwise, there is an intermediate path $L = yz$ such that $y \in V(I_1^*)$ and $z \in V(I_2^*)$. If $|L| \geq 2$ then by Lemma 2,

$$|C| = |I_1| + |I_2| \geq 2\overline{p} + 2|L| + 4 \geq 2\overline{p} + 8 = 3\delta - \overline{p} + 2 = 2\delta + 4,$$

contradicting (1). hence $|L| = 1$, implying that $\Upsilon(I_1, I_2) \subseteq E(G)$. If there are two independent intermediate edges between $I_1, I_2$, then by Lemma 2, $|C| = |I_1| + |I_2| \geq 2\overline{p} + 8 = 3\delta - \overline{p} + 2 = 2\delta + 4$, contradicting (1). Otherwise $\tau \leq 1$, contradicting the hypothesis.

**Case 5.** $2 \leq \overline{p} = \delta - 1$.

It follows that $|N_C(x_i)| \geq \delta - \overline{p} = 1$ $(i = 1, 2)$.

**Case 5.1.** $|N_C(x_i)| \geq 2$ $(i = 1, 2)$.

If $N_C(x_1) \neq N_C(x_2)$ then by Lemma 1, $|C| \geq 2\overline{p} + 8 = 3\delta - \overline{p} + 5 > 2\delta + 4$, contradicting (1). Hence, $N_C(x_1) = N_C(x_2)$. Clearly $s \geq 2$. Further, if $s \geq 3$ then

$$|C| \geq s(\overline{p} + 2) \geq 3(\delta + 1) > 2\delta + 4,$$

contradicting (1). Let $s = 2$. Since $\kappa \geq 3$, there is an edge $zw$ such that $z \in V(P)$ and $w \in V(C) \setminus \{\xi_1, \xi_2\}$. Assume w.l.o.g. that $w \in V(I_1^*)$. Then it is easy to see that $|I_1| \geq \delta + 3$. Since $|I_2| \geq \delta + 1$, we have $|C| \geq 2\delta + 4$, contradicting (1).

**Case 5.2**. Either $|N_C(x_1)| = 1$ or $|N_C(x_2)| = 1$.

Assume w.l.o.g. that $|N_C(x_1)| = 1$. Put $N_C(x_1) = \{y_1\}$. If $N_C(x_1) \neq N_C(x_2)$ then $x_2y_2 \in E(G)$ for some $y_2 \in V(C)\backslash\{y_1\}$ and we can argue as in Case 4.1. Let $N_C(x_1) = N_C(x_2) = \{y_1\}$. Since $\kappa \geq 1$, there is an edge $zw$ such that $z \in V(P)$ and $w \in V(C)\backslash\{y_1\}$. Clearly, $z \notin \{x_1, x_2\}$ and $x_2z^- \in E(G)$, where $z^-$ is the previous vertex of $z$ along $\overrightarrow{P}$. Then replacing $P$ with $x_1\overrightarrow{P}z^-x_2\overleftarrow{P}z$, we can argue as in Case 4.1.

**Case 6**. $\overline{p} \geq \delta$.

If $|C| \geq \kappa(\delta + 1)$ then clearly $|C| \geq 2\delta + 4$, contradicting (1). Otherwise, by Lemma 3, we can assume that $|N_C(x_i)| \geq 2$ $(i = 1, 2)$. Then $|C| \geq 2(\overline{p} + 2) \geq 2\delta + 4$, contradicting (1). ∎

Institute for Informatics and Automation Problems
National Academy of Sciences
P. Sevak 1, Yerevan 0014, Armenia
E-mail: zhora@ipia.sci.am